


\font\eightrm=cmr8

\font\eightit=cmti8 
\font\eightbfti=cmbxti10 at 8pt
\font\bfti=cmbxti10
\font\sc=cmcsc10
\font\eighttt=cmtt8

\font\tenbb=msbm10
\font\sevenbb=msbm7
\font\fivebb=msbm5
\newfam\bbfam
\textfont\bbfam=\tenbb
\scriptfont\bbfam=\sevenbb
\scriptscriptfont\bbfam=\fivebb
\def\bb{\fam\bbfam\tenbb}
\let\oldbb=\bb
\def\bb #1{{\oldbb #1}}

\font\fourteenbf=cmbx10 at 14pt
\font\eightbf=cmbx8
\font\sixbf=cmbx6

\def\scs{\scriptstyle}


\def\RR{{\bb R}}
\def\ZZ{{\bb Z}}
\def\QQ{{\bb Q}}
\def\AA{{\bb A}}
\def\Real{\mathop{\rm Re}}

\def\Pf{\mathop{\rm PF}}


\overfullrule=0pt
\vsize= 20.5cm
\hsize=15cm
\hoffset=3mm
\voffset=3mm

\headline{\hfil\tentt PREPUBLICATION 594 -- OCTOBRE 2000 -- UNIVERSITE DE
  NICE\hfil math/0101068}

\vskip .5cm

{\parindent=0pt
\eightbf

{\eightrm C. R. Acad. Sci. Paris, t.~331, S\'erie~I,
p. 423-428, 2000


Th\'eorie des Nombres/\hskip -.5mm{\eightbfti Number Theory}

(Analyse Harmonique/\hskip -.5mm{\eightbfti Harmonic Analysis})

\vskip 1.5cm


{\fourteenbf
Sur les Formules Explicites I: analyse invariante
}
\parindent=-1mm\footnote{}{
Note pr\'esent\'ee par Jean--Pierre K{\sixbf AHANE.}
}

\vskip 15pt


\noindent{\bf Jean--Fran\c{c}ois BURNOL}

\vskip 6pt


{\parindent=3mm
Laboratoire J.~A.~Dieudonn\'e, Math\'ematiques, Universit\'e de
Nice~Sophia-Antipolis, Parc Valrose, F--06108 Nice cedex~02, France

Courriel~: burnol@math.unice.fr
\par}

\medskip


(Re\c{c}u le 21 juillet 2000, accept\'e le 17 ao\^ut 2000)
}

\vskip3pt

\line{\hbox to 2cm{}\hrulefill\hbox to 1.5cm{}}

\vskip3pt


{\parindent=2cm\rightskip 1.5cm\baselineskip=9pt
\item{\bf R\'esum\'e.~}{\eightrm Weil a montr\'e que les formules explicites qui
relient les id\'eaux premiers d'un corps de nombres~$\scs K$ aux z\'eros et
p\^oles des s\'eries~$\scs L$ de Dirichlet--Hecke font intervenir les
compl\'etions~$\scs K_\nu$ de~$\scs K$. Nous montrons que l'analyse de Fourier
(multiplicative) de la transformation de Fourier (additive) permet d'obtenir ce
r\'esultat, en traitant identiquement chaque place, finie ou infinie, ramifi\'ee
ou non. Par ailleurs nous v\'erifions le crit\`ere de positivit\'e de Weil sous
une condition de support. ~\copyright~2000 Acad\'emie des Sciences/\'Editions
scientifiques et m\'edicales Elsevier SAS}
\par}

\vskip 10pt


{\parindent=2cm\rightskip 1.5cm
{\bfti On the Explicit Formulae I: invariant analysis}
\par}

\vskip 10pt


{\parindent=2cm\rightskip 1.5cm\baselineskip=9pt
\item{\bf Abstract.~}{\eightit Weil has generalized the Riemann--von Mangoldt
explicit formula linking the prime numbers with the zeros of the zeta function
to the set--up of a general algebraic number field~$\scs K$ and
Dirichlet--Hecke~$\scs L$--function, revealing in the process the r\^ole played
by the completions (finite and infinite) of~$\scs K$. We show how the local
terms of these explicit formulae are explained by the dilation invariant
``conductor operator''~$\scs \log(|x|_\nu) + \log(|y|_\nu)$. We also check
Weil's positivity criterion under a support
condition.
{\eightrm ~\copyright~2000 Acad\'emie des Sciences/\'Editions
scientifiques et m\'edicales Elsevier SAS}}
\par}

}

\vskip1.5pt
\line{\hbox to 2cm{}\hrulefill\hbox to 1.5cm{}}
\medskip

\parindent=3mm
\vskip 5mm

\noindent{\bf Abridged English Version}

\bigskip

Let $K$ be a number field and $\chi$ a Hecke character of the idele group
$\AA^\times$ with local component $\chi_\nu$ at the place $\nu$. Haar measures
$dx$ on $K_\nu$ and $d^*t$ on $K_\nu^\times$ are normalized as in [5]. The 
Hilbert spaces $L^2(K_\nu,dx)$ and $L^2(K_\nu^\times,d^*t)$ are canonically
isometric. We write $q_\nu$ for the cardinality of the residue field at a
finite place.
The dual $X_\nu$ of $K_\nu^\times$ is a collection of circles (or lines if $\nu$
is archimedian). We parametrize the component $X_\nu^\chi$ containing
$\chi_\nu^{-1}$ with the help of $D = \{\Real(s) = {1\over2}\}$, sending $s$ to
$\chi_\nu^{-1}(t)|t|_\nu^{-(s-{1\over2}) }$. Let $\cal F$ be the additive
Fourier transform on $L^2(K_\nu,dx)$ and $I$ the inversion $\varphi(x) \mapsto
{1\over|x|_\nu}\varphi({1\over x})$. The composite $\Gamma = {\cal F}I$ commutes
with $K_\nu^\times$, hence is given by a spectral multiplier
$\Gamma(s,\chi_\nu)$ on $X_\nu^\chi$. Clearly these are the Tate Gamma functions
occuring in ${\cal F}(\chi_\nu(x)|x|_\nu^{s-1}) = \Gamma(s,\chi_\nu)
\chi_\nu^{-1}(x)|x|_\nu^{ - s}$ (for $0<\Real(s)<1$). Let $g$ be a smooth
compactly supported function on $(0,\infty)$ with Mellin transform
$\widehat{g}(s) = \int_0^\infty g(u) u^{s-1} du$. 
Let $W(g;\chi)$ be the sum with multiplicities
of $\widehat{g}(s-{1\over2})$ over the zeros and poles of the complete $L$--
function $L(s,\chi)$ (poles being counted negatively).

\medskip

{\sc Lemma.~--~}$$W(g;\chi) = \sum_\nu \int_{s={1\over2} + i\tau} \widehat{g}
(i\tau){\Gamma^\prime(s,\chi_\nu)\over\Gamma(s,\chi_\nu)}\,{d\tau\over2\pi}$$

\smallskip

The contribution due to the discriminant of $K$ is distributed among the finite
places. When $\nu$ is finite the integrand is periodic except for
$\widehat{g}(i\tau)$. Using Poisson summation

$${1\over \log(q_\nu)}\,\sum_{j\in\ZZ}\widehat{g}(j\,{2\pi i\over\log(q_\nu)}
 + i\tau) = \sum_{k\in\ZZ}g(q_\nu^k)\,q_\nu^{k\,i\tau} = 
\int_{K_\nu^\times}g(|t|_\nu) |t|_\nu^{i\tau}\,d^*t$$

\noindent and this integral gives the spectral decomposition of
$g(|t|_\nu)\chi^{-1}_\nu(t)\in L^2(K_\nu^\times,d^*t)$ on the
characters $\chi_\nu^{-1}(t)|t|_\nu^{-i\tau}$ in $X_\nu^\chi$. Let $A =
\log(|x|_\nu)$ acting on $L^2(K_\nu,dx)$ as $\varphi(x)\mapsto
\log(|x|_\nu)\varphi(x)$ and $B = {\cal F}A{\cal F}^{-1} = \log(|y|_\nu)$. On
$L^2(X_\nu)$ one has $A = {1\over i}
{\partial\over\partial\tau} = {\partial\over\partial s}$. Then $\Gamma A
\Gamma^{-1}$ acts as $\alpha(s)\mapsto \Gamma(s,\chi_\nu){\partial\over\partial
s}\left(\Gamma(s,\chi_\nu)^{-1}\alpha(s)\right)= {\partial\over\partial
s}\alpha(s) - {\Gamma^\prime(s,\chi_\nu)\over\Gamma(s,\chi_\nu)}\alpha(s)$. So
the invariant operator $A - \Gamma A \Gamma^{-1} = A + B$ has spectral
multipliers  $+{\Gamma^\prime(s,\chi_\nu)\over\Gamma(s,\chi_\nu)}$. We conclude:

\bigskip

{\sc Theorem.~--~} {\it Let $g_{\nu\,;\,\chi}(t) =
g(|t|_\nu)\chi^{-1}_\nu(t)\in L^2(K_\nu^\times,d^*t)$. Then
$$W(g;\chi) = \sum_\nu H_\nu(g_{\nu\,;\,\chi})(1)$$
where the {\rm conductor operator} $H_\nu$ acts as $\log(|x|_\nu) + \log(|y|_\nu)$
 on $L^2(K_\nu,dx)$.}

\medskip

{\it Note.--~} As $H_\nu$ commutes with ${\cal F}$ and with $\Gamma$ it also
commutes with $I$ (which on $L^2(K_\nu^\times,d^*t)$ is
$f(t)\mapsto f({1\over t}))$. This is a local manifestation of the global
functional equation.

\medskip

The theorem explains Weil's discovery~[6] of a $\nu$--adic origin of
the local terms of the Explicit Formula. Weil's formulae
for the local terms as well as Haran's~[3] (in the case of the Riemann 
zeta function) are consequences. More details and
background are given in the french--language section and in [1].

\bigskip

Let $Z(g) = \sum_{\rho} \widehat{g}(\rho)$ be the sum over the critical zeros of
the Riemann zeta function. Let $g^\tau(u) = \overline{{1\over u}g({1\over u})}$
and $k = g*g^\tau$. Then $Z(k) = \sum_{\rho}
\widehat{g}(\rho)\overline{\widehat{g}(\overline{1-\rho})}$ and it is elementary
that the Riemann Hypothesis is equivalent to: $Z(k)\geq 0$ for all smooth
compactly supported $g$'s (Weil~[6], for a wider class of $g$'s).

\medskip

{\sc Theorem.~--~} {\it There is a $c>1$ such that $Z(k)\geq 0$ for all smooth
$g$'s with support in $[{1\over c}, c]$.}

\smallskip

{\it Proof.~--~} For $c\leq\sqrt{2}$ the support of $k$ is contained in
$[{1\over 2}, 2]$ and 

$$Z(k) = 2\Real(\widehat{k}(0)) + \int_{s={1\over2} + i\tau}
h_+(\tau)|\widehat{g}(s)|^2\,{d\tau\over2\pi}$$

\noindent with $h_+(\tau) = -\log(\pi) +
\Real(\lambda({1\over4}+{1\over2}i\tau))$ (and $\lambda(s) =
\Gamma^\prime(s)/\Gamma(s)$). As $\widehat{k}(0) = \int_{{1\over2}}^\infty
k(u)\,{du\over u} = \int_D \widehat{k}(s)\,2^s\,{1\over s}\,{d\tau\over2\pi}$
and $0 = \int_D \widehat{k}(s)\,2^{1-s}\,{1\over s}\,{d\tau\over2\pi}$, we have
$2\Real(\widehat{k}(0)) = \int_D {8\sqrt{2}\cos(\log(2)\tau)\over 1 +
4\tau^2}|\widehat{g}(s)|^2\,{d\tau\over2\pi}$.

$$Z(k) = \int_{s={1\over2} + i\tau} \left({8\sqrt{2}\cos(\log(2)\tau)\over 1 +
4\tau^2} + h_+(\tau)\right)|\widehat{g}(s)|^2\,{d\tau\over2\pi}$$

\noindent The continuous (real-valued) function $\alpha(\tau) =
{8\sqrt{2}\cos(\log(2)\tau)\over 1 + 4\tau^2} + h_+(\tau)$ satisfies
$\lim_{\tau\to\pm\infty} \alpha(\tau) = +\infty$ hence for all sufficiently
small $\varepsilon>0$ and suitable $A_\varepsilon>0$ : $\forall \tau\ 
A_\varepsilon\cos(\varepsilon\tau) + \alpha(\tau) \geq 0$. Now $\int
\cos(\varepsilon\tau)|\widehat{g}(s)|^2\,{d\tau\over2\pi} = \Real\left(\int
e^{-i\varepsilon\tau}\widehat{k}({1\over2} + i\tau)\,{d\tau\over2\pi}\right) =
\Real(e^{\varepsilon/2}k(e^{\varepsilon}))$. For $c = e^{\varepsilon/2}$ one
then has $\int \cos(\varepsilon\tau)|\widehat{g}(s)|^2\,{d\tau\over2\pi} = 0$
for $g$ with support in $[{1\over c}, c]$, hence $Z(k)\geq 0$. Computer
calculations help being more precise about the allowable $c$'s but anyhow a
further idea seems necessary to reach $c = \sqrt{2}$.


\vskip 6pt
\centerline{\hbox to 2cm{\hrulefill}}
\vskip 9pt



\noindent{\bf 1. Op\'erateurs invariants}

\medskip

Il est bien connu que tout op\'erateur born\'e sur $L^2(\RR,dx)$ qui commute
avec les translations est diagonalis\'e par la transformation de Fourier mais il
nous est utile d'autoriser des multiplicateurs non born\'es (et des groupes plus
g\'en\'eraux). Nous regroupons ici quelques lemmes dans ce sens \`a d\'efaut de
conna\^{\i}tre une r\'ef\'erence adapt\'ee (voir [2] pour les d\'emonstrations).
Soient $G$ un groupe topologique s\'epar\'e, localement compact et ab\'elien (voir
Rudin~[4]) et $\widehat{G}$ son dual. Il existe sur $G$ une mesure de Haar $dx$
et sur $\widehat{G}$ la mesure duale $dy$ pour lesquelles la transformation
$F(\varphi)(y) = \int \varphi(x) \overline{y(x)} dx$ est une isom\'etrie de
$L^2(G, dx)$ sur $L^2(\widehat{G}, dy)$. Nous supposerons que $dy$ est une
mesure $\sigma$-finie. Pour $\varphi \in L^2(G, dx)$ et $g \in G$ on note
$g\cdot\varphi$ la fonction $x\mapsto \varphi(x - g)$. Un op\'erateur $M$ de
domaine $D\subset L^2(G, dx)$ {\it commute avec $G$} si $\forall g\ \forall
\varphi: \varphi\in D\Rightarrow \left(\ g\cdot\varphi\in D\ \hbox{et}\
M(g\cdot\varphi) = g\cdot M(\varphi)\ \right)$. Soit $a(y)$ mesurable (finie
presque partout), $D_a = \{\varphi \in L^2(G, dx)\,|\,a\cdot F(\varphi)\in
L^2(\widehat{G}, dy)\}$, et $M_a$ de domaine $D_a$ agissant selon $M_a(\varphi)
= F^{-1}(a\cdot F(\varphi))$. On consid\`ere que $a=b$ si $a(y) = b(y)\,p{.}p$.

\medskip

{\sc Lemme~1.1.~--~}{\it Le domaine $D_a$ est dense dans $L^2(G, dx)$,
l'op\'erateur $(M_a, D_a)$ est ferm\'e et commute avec $G$. Tout op\'erateur de
domaine dense dans $L^2(G, dx)$, qui commute avec $G$, et qui est ferm\'e, est
(uniquement) de la forme $(M_a, D_a)$.}

\medskip

{\sc Corollaire~1.2.~--~}{\it Soit $(M,D)$ un op\'erateur de domaine $D$ dense
dans $L^2(G, dx)$, qui commute avec $G$, et qui est sym\'etrique: $\forall
\varphi,\psi\in D \ \int \overline{\varphi(x)}\,M(\psi)(x) dx =
\int\overline{M(\varphi)(x)}\,\psi(x) dx$. Alors $(M,D)$ est essentiellement
auto-adjoint, et sa cl\^oture auto-adjointe est l'unique $(M_a, D_a)$
v\'erifiant $(M_a, D_a)\supset (M,D)$.}

\medskip

{\sc Lemme~1.3.~--~}{\it Soit $L$ un espace de Hilbert et $G$ un groupe
d'op\'erateurs unitaires sur $L$ (non n\'ecessairement ab\'elien). Soit $M$ un
op\'erateur de domaine $D$ dense dans $L$, sym\'etrique, et commutant avec $G$.
Si l'alg\`ebre de von Neumann des op\'erateurs born\'es qui commutent avec $G$
est ab\'elienne alors $(M,D)$ est essentiellement auto-adjoint.}

\bigskip

\noindent{\bf 2. Formules de Weil}

\medskip

Nous exposons certains points du travail de Weil~[6], pour une s\'erie
$L(s,\chi)$ de Dirichlet. Notons $\rho$ les z\'eros (avec multiplicit\'es) de
$L(s,\chi)$ avec $0\leq\Real(s)\leq1$. Soit $\nu$ une place de $\QQ$, donc $\nu =
p$ (correspondant \`a $\QQ_p$) ou $\nu = r$ correspondant \`a $\RR$. Soit $\AA$
l'anneau des ad\`eles et $\AA^\times$ le groupe des id\`eles. Notons
$U_p\subset\QQ_p^\times$ le sous-groupe des unit\'es. On a $\AA^\times \cong
\QQ^\times\times\RR^{\times+}\times\prod_p U_p$ o\`u chaque terme de droite est
identifi\'e \`a son plongement dans $\AA^\times$. Ainsi un caract\`ere de Hecke
$\chi_H:\ \AA^\times \rightarrow U(1)$ (continu, trivial sur $\QQ^\times$), est
la donn\'ee d'un caract\`ere $u\mapsto u^{i\tau}$ de $\RR^{\times+}$, et d'un
nombre fini de caract\`eres non triviaux des $U_p$. Ceux-ci \'equivalent \`a la
donn\'ee d'un caract\`ere $\chi$ {\it primitif} de $(\ZZ/q\ZZ)^\times$ (pour $q
= 1$ lire $\{1\}$) pour un certain conducteur $q$. Notons $\pi_p$ l'id\`ele de
composantes $1$ pour $\nu\neq p$ et $p$ pour $\nu = p$. On aura
$p\not{}\hskip-2pt |\ q\Rightarrow \chi_H(\pi_p) = \chi(p)^{-1}p^{-i\tau}$. On
associe donc au caract\`ere de Dirichlet $\chi$ un caract\`ere de Hecke, encore
not\'e $\chi$, qui sera $\chi_H^{-1}$ (avec $\tau = 0$). Soit $\chi_\nu$ la
composante locale obtenue par $\QQ_\nu^\times \hookrightarrow \AA^\times
{\buildrel\chi\over\rightarrow}\, U(1)$. Alors $\chi_p(p) = \chi(p)$ lorsque $p$
est premier avec $q$. Par ailleurs $\chi_r(t) = 1$ pour $t>0$, $\chi_r(t) =
\chi(-1)$ pour $t<0$. Soit $g(u)$ une fonction \`a support compact sur
$(0,\infty)$, de classe ${\cal C}^\infty$. Notons $\widehat{g}(s)$ sa
transform\'ee de Mellin $\int_0^\infty g(u) u^{s-1} du$. Posons $\delta_\chi =
1$ si $\chi$ est le caract\`ere principal, $ = 0$ sinon.

\medskip

{\sc Th\'eor\`eme~2.1.~(Weil~[6])--~}{$$\sum_\rho \widehat{g}(\rho-{1\over2}) -
\delta_\chi\,\left(\,\widehat{g}(-{1\over2})+\widehat{g}({1\over2})\right)
 = - \sum_\nu
\Pf{}_\nu\int_{\QQ_\nu^\times} {g({1\over |t|_\nu})\sqrt{|t|_\nu}\,\chi_\nu(t)
\over |1 -t|_\nu}\,d^\times t$$}

\smallskip

Dans cette formule $d^\times t$ est la mesure de Haar sur $\QQ_p^\times$ qui
donne une masse de $\log(p)$ \`a $U_p$. Pour $\RR^\times$ il s'agit de
${dt\over2|t|}$. La partie finie $\Pf_\nu$ d\'esigne une certaine
r\'egularisation en $t = 1$. La somme sur les z\'eros (et les p\^oles)
s'obtient par une int\'egrale de contour sur des rectangles et fait donc
intervenir la d\'eriv\'ee logarithmique de $\xi(s,\chi) =
\pi^{-s/2}\Gamma({s\over2})L(s,\chi)$ (nous supposerons  que $\chi$ est {\it
pair\/}: $\chi(-1) = 1$). Par transformation de Mellin inverse on obtient une
somme $\sum_\nu w_\nu(g;\chi)$ o\`u les $w_\nu(\ ;\chi)$ sont {\it a priori} des
distributions sur $(0,\infty)$. Pour la place archim\'edienne le bord droit du
rectangle donne $-g(1)(\log(\pi)+\gamma)/2 - \int_1^\infty g(u)u^{-{1\over2}}
{du\over u} -
\int_1^\infty {g(u)u^{-{1\over2}} - g(1) \over u^2 - 1}{du\over u}$.
Avec la contribution du
bord gauche cela donne myst\'erieusement (et la formule vaut aussi pour
$\chi(-1) = -1$):

$$w_r(g;\chi) = - (\log(2\pi)+\gamma)g(1) - \int_{1\over2}^\infty {g({1\over t})
\sqrt{t}
- g(1)\over |1 - t|}{d^\times t} - \int_0^{1\over2} {g({1\over t})\sqrt{t}
 \over |1 -t|}{d^\times t} - \int_{-\infty}^{0}{g({1\over|t|})\sqrt{|t|}
\chi(-1)\over |1 - t|}{d^\times t}\leqno\hbox{\sc (2.2.)} $$

\noindent Pour $p\not{}\hskip-1pt |\ q$ on obtient $- \log(p)\sum_{j\geq 1}
\left(\chi(p)^j {p}^{-{j\over2}}g(p^j)+ \chi(p)^{-j}{p}^{-{j\over2}}g(p^{-j})\right)$
 et dans le cas
ramifi\'e $p\,|\,q$ on obtient par l'\'equation fonctionnelle $w_p(g;\chi) =
+f_p(\chi) \log(p)\cdot g(1)$ o\`u $f_p(\chi)$ est l'exposant de $p$ dans $q$.
Or un calcul explicite \'el\'ementaire donne:

\medskip

{\sc Th\'eor\`eme~2.3.~(Weil~[6])--~}{$$ f_p(\chi) \log(p) =  \int_{|t|_p = 1}
{1 - \chi_p(t)\over |1 - t|_p}\,d^\times t$$}

\smallskip

\noindent ce qui permet \`a Weil d'\'ecrire d'une mani\`ere g\'en\'erale pour
une place finie

$$w_p(g;\chi) =  - \int_{|t|_p\neq1} {g({1\over |t|_p})\sqrt{|t|_p}\chi_p(t)
\over |1 - t|_p}\,d^\times t - \int_{|t|_p = 1} {g(1)\chi_p(t) - g(1)\over |1 -
t|_p}\,d^\times t\leqno\hbox{\sc (2.4.)}$$

Haran~[3] a montr\'e (pour la fonction $\zeta(s)$ de Riemann) que l'on pouvait
reformuler les $w_\nu(g;\chi)$ comme des convolutions {\it additives}. Soit
$R_s^\nu$ le noyau de convolution additive sur $\QQ_\nu$ dont la transformation
de Fourier est $|y|_\nu^{-s}$. 

\medskip

{\sc Th\'eor\`eme~2.5.~(Haran~[3]{\rm\ pour $\chi=\hbox{\bf 1}$})--~}{\it Posons
$\varphi_{\nu\,;\,\chi}(x) = g(|x|_\nu)|x|_\nu^{-{1\over2}}\chi^{-1}_\nu(x)$ pour
 $x\in\QQ_\nu$, $x\neq0$ et $\varphi_{\nu\,;\,\chi}(0) = 0$. Alors
 $w_\nu(g;\chi) = -\left.{\partial\over\partial s}\right|_{s=0}
 (R_s^\nu*\varphi_{\nu\,;\,\chi})(1)$.}

\smallskip

Nous mettrons moins en avant le ``semi--groupe de Riesz'' $R_s^\nu$ que Haran,
au profit de la transform\'ee de Fourier $G_\nu(x)$ de la distribution
$\log(|y|_\nu)$. En effet {\sc (2.5.)} \'equivaut \`a:

$$w_\nu(g;\chi) = + (G_\nu * \varphi_{\nu\,;\,\chi})(1)\leqno\hbox {\sc (2.6.)}$$

\noindent Il est important que {\sc (2.6.)}\ soit \'egalement valable pour $q>1$
et $\nu = p$, $p\,|\,q$ (cas ramifi\'e). On peut calculer $G_\nu(x)$ en
utilisant ${\cal F}(|x|_\nu^{s-1}) = \Gamma(s,\hbox{\bf 1})|x|_\nu^{ - s}$ pour
$s=1-\varepsilon, \varepsilon\to0$ (avec $\Gamma(s,\hbox{\bf 1})$ la fonction
Gamma de Tate~[5]). On v\'erifie alors que {\sc (2.6.)} est identique avec {\sc
(2.2.)} et {\sc (2.4.)}.

\bigskip

\noindent{\bf 3. L'op\'erateur conducteur}

\medskip

Nous obtenons {\sc (2.6.)} simultan\'ement pour les places r\'eelles, complexes,
finies et non-ramifi\'es, finies et ramifi\'ees (voir [1] pour un expos\'e plus
d\'etaill\'e et des r\'ef\'erences suppl\'ementaires). Il serait int\'eressant
d'\'etendre notre interpr\'etation de {\sc (2.3.)} aux caract\`eres
non-ab\'eliens de Artin ([7], [8]).

Soit $K$ un corps de nombres et $\chi$ un caract\`ere de Hecke du groupe des
id\`eles $\AA^\times$ de $K$. Soit $\nu$ une place de $K$. Si $\nu$ est finie 
on note $q_\nu$ la cardinalit\'e du corps r\'esiduel. Soit $dx$ la mesure
de Haar sur $K_\nu$, normalis\'ee comme dans la th\`ese de Tate~[5]. Soit $d^*t$
la mesure de Haar sur $K_\nu^\times$, pour laquelle les unit\'es ont volume $1$
(place finie) ou qui est ${dt\over2|t|}$ (place r\'eelle) ou encore
$drd\theta\over\pi r$ (place complexe). Les espaces $L^2(K_\nu,dx)$ et
$L^2(K_\nu^\times,d^*t)$ sont isom\'etriques. Le dual $X_\nu$ de $K_\nu^\times$
est une collection de cercles (ou de droites dans le cas archim\'edien)
index\'es par les caract\`eres du sous-groupe des unit\'es. Nous param\'etrons
la composante $X_\nu^\chi$ contenant $\chi_\nu^{-1}$ par la droite critique $D =
\{\Real(s) = {1\over2}\}$ en associant \`a $s$ le caract\`ere
$\chi_\nu^{-1}(t)|t|_\nu^{-(s-{1\over2}) }$. Soit $\cal F$ la transformation de
Fourier additive qui est donc une isom\'etrie de $L^2(K_\nu,dx)$, et $I$
l'inversion $\varphi(x) \mapsto {1\over|x|_\nu}\varphi({1\over x})$. Le
compos\'e $\Gamma = {\cal F}I$ commute avec l'action de $K_\nu^\times$, il lui
correspond donc des multiplicateurs spectraux $\Gamma(s,\chi_\nu)$. On v\'erifie
ais\'ement qu'il s'agit l\`a des fonctions Gamma de Tate qui apparaissent dans
l'identit\'e de distributions sur $K_\nu$: ${\cal F}(\chi_\nu(x)|x|_\nu^{s-1}) =
\Gamma(s,\chi_\nu) \chi_\nu^{-1}(x)|x|_\nu^{ - s}$ (pour $0<\Real(s)<1$). Soit
$g(u)$ une fonction \`a support compact sur $(0,\infty)$, de classe ${\cal
C}^\infty$. Notons $W(g;\chi)$ la somme avec multiplicit\'es de 
$\widehat{g}(s-{1\over2})$
sur les z\'eros et p\^oles de la fonction $L$ compl\`ete $L(s,\chi)$ 
(les p\^oles \'etant compt\'es n\'egativement).

\medskip

{\sc Lemme~3.1.~--~}$$W(g;\chi) = \sum_\nu \int_{s={1\over2} + i\tau}
\widehat{g}(i\tau)
{\Gamma^\prime(s,\chi_\nu)\over\Gamma(s,\chi_\nu)}\,{d\tau\over2\pi}$$

\smallskip

Dans cette expression la contribution due au discriminant de $K$ est ventil\'ee
sur les places finies. Lorsque $\nu$ est finie l'int\'egrand est p\'eriodique
\`a l'exception de $\widehat{g}(i\tau)$. Par sommation de Poisson

$${1\over \log(q_\nu)}\,\sum_{j\in\ZZ}\widehat{g}(j\,{2\pi i\over\log(q_\nu)}
 + i\tau) = \sum_{k\in\ZZ}g(q_\nu^k)\,q_\nu^{k\,i\tau} = 
\int_{K_\nu^\times}g(|t|_\nu) |t|_\nu^{i\tau}\,d^*t$$

\noindent On notera que la mesure sur $X_\nu^\chi$ duale \`a $d^*t$ est
$\log(q_\nu)\,{d\tau\over2\pi}$ (sur une p\'eriode de longueur
${2\pi\over\log(q_\nu)}$) et que la derni\`ere
int\'egrale donne la d\'ecomposition spectrale de $g(|t|_\nu)
\chi^{-1}_\nu(t)\in L^2(K_\nu^\times,d^*t)$ par les caract\`eres
$\chi_\nu^{-1}(t)|t|_\nu^{-i\tau}$ dans $X_\nu^\chi$.

\medskip

{\sc Lemme~3.2.~--~} {\it Soit $g_{\nu\,;\,\chi}(t) =
g(|t|_\nu)\chi^{-1}_\nu(t) \in L^2(K_\nu^\times,d^*t)$
et $H_\nu$ l'op\'erateur sur 
$L^2(K_\nu^\times,d^*t)$ de
multiplicateurs spectraux ${\Gamma^\prime(s,\chi_\nu)\over\Gamma(s,\chi_\nu)}$
sur les $X_\nu^\chi$. Alors $W(g;\chi) = \sum_\nu H_\nu(g_{\nu\,;\,\chi})(1)$}

\smallskip

Soit $A = \log(|x|_\nu)$ l'op\'erateur (non-born\'e) sur $L^2(K_\nu,dx)$ qui
agit selon $\varphi(x)\mapsto \log(|x|_\nu)\varphi(x)$ et soit $B = {\cal
F}A{\cal F}^{-1}$ l'op\'erateur conjugu\'e. On \'ecrira:  $B = \log(|y|_\nu)$.
Sur $L^2(X_\nu)$ on a $A_\nu = {1\over i}
{\partial\over\partial\tau}$, plus bri\`evement $A = {\partial\over\partial s}$.
L'op\'erateur $\Gamma A \Gamma^{-1}$ agit selon $\alpha(s)\mapsto
\Gamma(s,\chi_\nu){\partial\over\partial
s}\left(\Gamma(s,\chi_\nu)^{-1}\alpha(s)\right)= {\partial\over\partial
s}\alpha(s) - {\Gamma^\prime(s,\chi_\nu)\over\Gamma(s,\chi_\nu)}\alpha(s)$.
Ainsi $H_\nu = A - \Gamma A \Gamma^{-1} = A + B$. En conclusion:

\bigskip

{\sc Th\'eor\`eme~3.3.~([1])--~} {\it Soit $g_{\nu\,;\,\chi}(t) =
g(|t|_\nu)\chi^{-1}_\nu(t) \in L^2(K_\nu^\times,d^*t)$. On a
$$W(g;\chi) = \sum_\nu H_\nu(g_{\nu\,;\,\chi})(1)$$
o\`u $H_\nu$ est {\rm l'op\'erateur conducteur} qui s'\'ecrit
$$H_\nu = \log(|x|_\nu) + \log(|y|_\nu)$$
sur $L^2(K_\nu,dx)$.}

\bigskip

{\it Note~3.4.~--~} On peut \'evaluer $H_\nu(g_{\nu\,;\,\chi})(t)$ en $1$ sans
ambigu\"{\i}t\'e car elle admet un repr\'esentant lisse (${\cal C}^\infty$ pour
$\nu$ archim\'edienne, localement constant pour $\nu$ finie).

\medskip

{\it Note~3.5.~--~} Soit $\varphi_{\nu\,;\,\chi}\in L^2(K_\nu,dx)$ d\'efini
par $\varphi_{\nu\,;\,\chi}(x) = g(|x|_\nu)\chi^{-1}_\nu(x)$ pour $x\neq0$
(et $\varphi_{\nu\,;\,\chi}(0) =0$). 
Consid\'erons $H_\nu$ comme un op\'erateur sur
$L^2(K_\nu,dx)$. Alors $\sum_\nu H_\nu(\varphi_{\nu\,;\,\chi})(1)$ est la
somme de $\widehat{g}(s)$ sur les z\'eros et p\^oles de $L(s,\chi)$.

\medskip

{\it Note~3.6.~--~} $H_\nu$ commute avec $I$: en effet il commute avec $\Gamma$ et
avec ${\cal F}$. On obtient l\`a une manifestation locale de l'\'equation
fonctionnelle globale des s\'eries $L$.

\medskip

Revenons au corps $\QQ$ et notons $Z(g) = \sum_{\rho} \widehat{g}(\rho)$ la
somme sur les z\'eros (non-triviaux) de $\zeta(s)$ . Soit $g^\tau(u) =
\overline{{1\over u}g({1\over u})}$ de sorte que $\widehat{g^\tau}(s) =
\overline{\widehat{g}(\overline{1-s})}$. Soit $k = g*g^\tau$ la convolution
multiplicative. Alors $Z(k) = \sum_{\rho}
\widehat{g}(\rho)\overline{\widehat{g}(\overline{1-\rho})}$ et il est
\'el\'ementaire que l'Hypoth\`ese de Riemann \'equivaut \`a: $Z(k)\geq0$ pour
toute fonction $g$ de classe ${\cal C}^\infty$ \`a support compact (Weil~[6],
pour une classe plus large de fonctions $g$).

\medskip

{\sc Th\'eor\`eme~3.7.~--~} {\it Il existe $c>1$ tel que $Z(k)\geq 0$ pour toute
fonction $g$ de classe ${\cal C}^\infty$ \`a support dans $[{1\over c}, c]$.}

\smallskip

La d\'emonstration est donn\'ee dans la partie en langue anglaise de cette Note.

\bigskip

\noindent{\bf Remerciements.} 
Je remercie la ``Soci\'et\'e de secours des amis
des sciences'' (SSAS) pour l'aide qu'elle m'a apport\'ee en 1999.

\bigskip


\centerline{\bf R\'ef\'erences bibliographiques}

\medskip

{\baselineskip=9pt
\eightrm
\item{[1]}Burnol~J.-F., math/9902080. 
En pr\'eparation. Une version pr\'eliminaire,
intitul\'ee ``The Explicit Formula and the conductor operator'' (f\'evrier 1999)
28 pp{.}, est disponible \`a l'adresse {\eighttt http://arXiv.org/abs/math/9902080}

\smallskip

\item{[2]}Burnol~J.-F., math/9907013. 
Addendum to ``Quaternionic gamma functions'' (juillet
1999) 6 pp{.}, manuscrit \'electronique disponible \`a l'adresse
{\eighttt http://arXiv.org/abs/math/9907013}

\smallskip

\item{[3]}Haran~S., Riesz potentials and explicit sums in arithmetic, Invent.
Math.  101 (1990) 697--703.

\smallskip

\item{[4]}Rudin~W., Fourier analysis on groups, Interscience Publishers, 1962.

\smallskip

\item{[5]}Tate~J., ``Fourier analysis in number fields and Hecke's
zeta-function'' (Princeton 1950) dans Algebraic Number Theory, Proc.
Instructional Conf. Brighton, Academic Press, 1967, 305--347.

\smallskip

\item{[6]}Weil~A., Sur les ``formules explicites'' de la th\'eorie des nombres
premiers, Comm. S\'em. Math. Univ. Lund, volume d\'edi\'e \`a Marcel Riesz
(1952) 252--265. Oeuvres, Vol{.}~II. 

\smallskip

\item{[7]}Weil~A., Sur les formules explicites de la th\'eorie des nombres, Izv.
Mat. Nauk. (Ser. Mat.) 36 (1972) 3--18. Oeuvres, Vol{.}~III.

\smallskip

\item{[8]}Weil~A., Basic Number Theory, 3${}^{\rm i\grave eme}$ \'ed{.},
Springer Verlag New York, 1974.

}

\end